\documentclass[article,oneside,12pt]{article}

\input{epsf}
\usepackage{amsmath}
\usepackage{amsthm}
\usepackage{amsfonts}
\usepackage{amssymb}
\usepackage{mathrsfs}
\usepackage[dvips]{graphics}
\usepackage{graphicx}
\usepackage{float}
\usepackage{anysize}
\newtheorem{theorem}{Theorem}

\newcommand\tdots{,\dots ,}
\newcommand\rstr{R^{*}}

\newcommand{\xb}{\textbf{\textit{x}}}
\newcommand{\R}{\mathbb{R}}

\begin{document}

\title {Simulation Studies of Some Voronoi Point Processes%
\footnote{This research was supported by the ARC Centre of Excellence for Mathematics
and Statistics of Complex Systems.}}

\author{K.A. Borovkov\footnote{University of Melbourne, email: k.borovkov@ms.unimelb.edu.au}
and D.A. Odell\footnote{MASCOS, Melbourne, email: d.odell@ms.unimelb.edu.au}}

\maketitle

\begin{abstract}
We introduce a new class of dynamic point process models with simple and intuitive
dynamics that are based on the Voronoi tessellations generated by the processes. Under
broad conditions, these processes prove to be ergodic and produce, on stabilisation, a
wide range of clustering patterns. In the paper, we present results of simulation
studies of three statistical measures (Thiel's redundancy, van Lieshout and Baddeley's
$J$-function and the empirical distribution of the Voronoi nearest neighbours' numbers)
for inference on these models from the clustering behaviour in the  stationary regime.
In particular, we make comparisons with the area-interaction processes of Baddeley and
van Lieshout.
\end{abstract}

\section{Introduction}

In many diverse fields such as biology, geography, business and telecommunications,
spatial point configurations evolve according to criteria dependent on a {\em zone of
influence\/} in some general sense. A natural way of formalising this concept is based
on using Voronoi tessellations~(see e.g.~\cite{OBS}).

Let $(M,d)$ be a metric space and $N>0$ the initial number of points in our process.
The associated configuration space $\mathscr{X}$ consists of all $N$-point subsets of
$M$: $\mathscr{X}:=\{\xb\subset M :\, \text{card} (\xb)=N \}$ (`multiple points' will
be a.s. impossible in our models). For an $\xb\in \mathscr{X}$, the {\em Voronoi
cell\/} of $x_i$ relative to $\xb$ is defined as
$$
C^{\xb}_{x_i}:=  \Bigl\{y \in M: \, d(y,x_i) = \min_{j\le N}  d(y,x_j) \Bigr\}.
$$
The set $\mathscr{T}_{\xb}= \{C^{\xb}_{x_i}: \, x_i\in \xb \}$ is called the {\em
Voronoi tessellation\/} generated by $\xb$.

We define a Voronoi point process as a discrete-time Markov processes $\{\xb_n
\}_{n\geq 0}$ with values in $\mathscr{X}$ which evolves as follows: at each step, a
point is chosen from the current configuration $\xb_n$ according to a probability rule
determined by $\mathscr{T}_{\xb_n}$ and removed from the configuration, and at the same
time a new point is added at a random location, according to a fixed probability
measure $\mu$ on $M$. Our initial interest interest in such dynamics was prompted by
its relevance to some of the real-life processes in the above-mentioned application
areas. In addition, simulations showed that processes constructed according to this
model display very interesting forms of clustering behaviour, and constructing models
of point processes producing desired clustering patterns is of substantial interest both
theoretically and practically.~(see e.g.~\cite{BV95, OBS}).

Despite the fact that it is notoriously difficult to obtain analytic results for
processes of such a nature, we were able to prove ergodicity in a number of interesting
cases. To formulate the results, we first need to further specify the `culling rule' in
the model.

We assume that a non-negative `selection function'  $S$ is given on the set of all
possible Voronoi cells in $M$. Then, given that $\xb=\{x_{1},\dots ,x_{N}\}$ is the
current configuration of the point process, in the next step a random point $x_{J}\in
\xb$ is chosen to be removed, with the distribution of the index $J$ given by
\begin{equation}
\label{eq2}
 \mathbf{P} (J=j|\xb)
 :=\frac{S(C_{x_{j}}^{\xb})}{\sum_{i=1}^{N}S(C_{x_{i}}^{\xb})},
  \qquad  1\leq J \leq N.
\end{equation}
The function $S$ can be  based on different properties of a Voronoi cell. In this note,
we consider only two of them: the `volume' of the cell and the number of its edges.
More precisely, we introduce the following two classes of Voronoi point processes:

\medskip

(A) The volume-based Voronoi point process, or $v$-process, on  $M=\mathbb{S}^1$ or $M=
[0,1]^2$. Let $\lambda$ be the `volume measure' on $M$ (length on $\mathbb{S}^1$, area
on $[0,1]^2$). We assume that the value of the selection function $S(C_{x_{j}}^{\xb})$
is determined by the volume of the cell $C_{x_{j}}^{\xb}$: for a function
$S_v:\R^{+}\mapsto\R^{+}$, one puts $S(C_{x_{j}}^{\xb}):=S_{v}(\lambda
(C_{x_{j}}^{\xb}))$. If $S_v$ is increasing, then points with Voronoi cells of large
volume are more likely  to be culled, and so the selection favours points with small
cells, \text{i.e.} points restricted by `close neighbours'. A decreasing $S_v$ favours
points with large cells. Functions of the form $S_v(u)=u^{\alpha}$, $\alpha \in \R$,
produce scale-independent models (in this case the dynamics of the process will clearly
be invariant under scale transformations of~$M$).

\medskip

(B) The neighbour-based Voronoi point process, or $n$-process, on $M= [0,1]^2$. The set
$\xb[x_{j}]$ of Voronoi-nearest neighbours of a point $x_j\in \xb$ is defined as the
collection of those generators $x_k\in \xb$ whose Voronoi cells share an edge with
$C_{x_{j}}^{\xb}$:
$$
\xb[x_{j}] :=\Bigl\{x_{k}\in \xb : \, k\neq j,\, {\rm card}(C_{x_{k}}^{\xb} \cap
C_{x_{j}}^{\xb})>1 \Bigr\}.
$$
We assume that, for a function $S_n:\{1\tdots N\}\mapsto\R^{+}$, we have
$S(C_{x_{j}}^{\xb})=S_{n}({\rm card}( \xb[x_{j}] ))$.

\medskip
In both cases (A) and (B), we assume the placement probability $\mu=\lambda$.

\begin{theorem}
{\em (i)} The neighbour-based Voronoi point process with a positive selection function
$S_n: \{1,\dots , N \}\mapsto \R^+$ is Harris ergodic.

{\em (ii)} The volume-based Voronoi point process with a selection function
$S_v:(0,\lambda(M)]\mapsto\R^+$, such that both $S_v$ and $1/S_v$ are bounded on closed
subsets, is Harris ergodic.
\end{theorem}

Recall that Harris ergodicity entails convergence to the stationary distribution in
total variation (see e.g.  p.560 in~\cite{MeTw92} or p.154 in~\cite{As87}). For the
proof of the theorem, see~\cite{borode}.

This result justifies the approach taken in the present note, which is devoted to
simulation studies of point patterns emerging in stationary regime in the evolution of
Voronoi point processes. More specifically, we are interested in the performance of
three relatively simple statistics one could employ to characterise the resulting
stationary distributions. These are Thiel's redundancy measure (which is essentially
the relative entropy for the set of cell volumes) and, in the two-dimensional case, the
distribution of the number of Voronoi nearest neighbours for a random cell  and the
$J$-function of Baddeley and van Lieshout (which provides a comparison of the
environment of a `random point' of the configuration with that of a `random point' of
the underlying space over a range of scales). We demonstrate  that the $n$- and
$v$-processes produce (for different choices of parameters) quite different clustering
patterns, and that the combination of the above three measures work reasonably well in
distinguishing between different Voronoi processes, and also between Voronoi processes
and the area-interaction processes of Baddeley and Van Lieshout~\cite{BV95}.

\section{Simulation Studies}
\label{s4}
\subsection{One-dimensional $v$-processes}

\begin{figure}[htbp]
    \centering
        \includegraphics{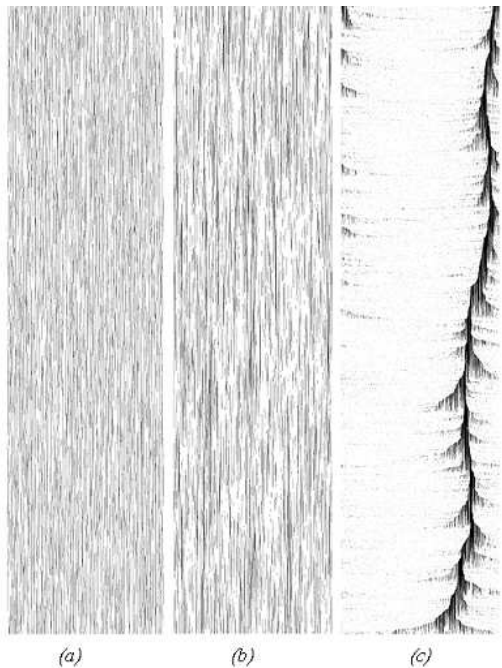}
        \caption{\footnotesize Evolution of $v$-processes on circle. Each of the three processes
        (with the selection functions $S_v(u)=u^{\alpha}$ with (a) $\alpha = -1.0$, (b)  $\alpha = 0.5$, (c)  $\alpha=1.5,$ resp.)
        was run with $N = 128$ points for $T = 4096$ steps (the vertical axis represents time). }
    \label{3oneD}
\end{figure}

The  simple model of  scale-free process on $M=\mathbb{S}^1$ with $S_v(u)=u^{\alpha}$
demonstrates that very interesting dynamics and point patterns arise even in one-dimension.
This includes a `phase change' observed when varying the value of
$\alpha$. Figure~\ref{3oneD} depicts side-by-side realisations of three different
$v$-processes of this type, each having the same number of points $N=128$, but
different $\alpha$ values. The base of each rectangle represents the circle
$\mathbb{S}^1$ opened out into a line segment by a cut.

We see a dramatic phase shift in behaviour between (b) and (c) that can be located more
precisely at $\alpha=1$; its sharpness was found to increase with~$N$. When $\alpha
>1$,  a cluster forms whose stability also increases with $N$. When $\alpha \leq
1$, we observe a degree of clustering in $M$ that varies with $\alpha$.

When $\alpha=0$, the points are uniformly distributed, so the associated Voronoi cell
volumes have pretty nearly a Gamma distribution with shape parameter 2 (as cell-width
is half the sum of two consecutive uniform spacings). For $\alpha \leq 1$, the
histogram of cell volumes stabilises after a large number of steps to a result
well-fitted by a Gamma distribution with a shape parameter varying inversely with
$\alpha$. This suggests that the maximum likelihood  estimator for the Gamma shape
parameter for the cell-volume data could be a  suitable statistic for inference on
$\alpha$. Furthermore, under the assumption of Gamma distribution, the shape parameter
can equivalently be estimated from the entropy---and the latter can also be used
without the above assumption, namely in the form of {\em Thiel's redundancy measure}
which is introduced as follows.

\begin{figure}[h]
    \centering
        \includegraphics{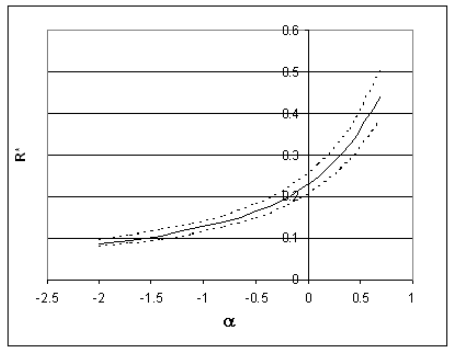}
        \caption{\footnotesize Empirical Thiel's redundancy measure for $v$-processes on
        $\mathbb{S}^1$ with
        $S_v(u)=u^{\alpha}$, $N=10^3$.
        The dotted lines show the 95\%\ probability intervals.}
    \label{Rstar1D}
\end{figure}

\begin{figure}[!ht]
    \centering
        \includegraphics{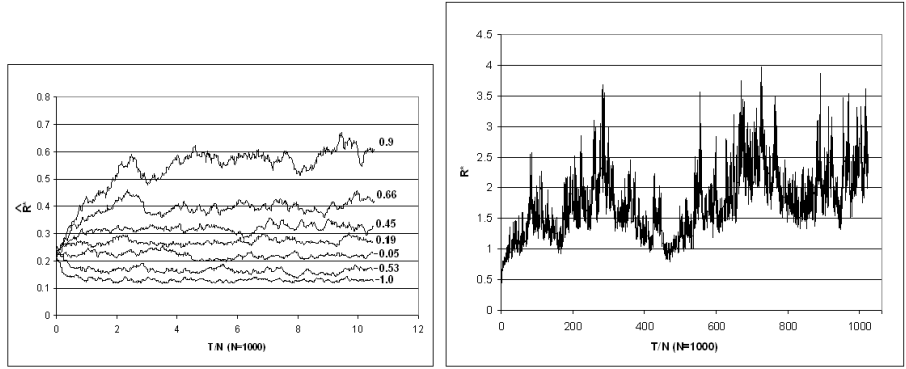}
        \caption{\footnotesize Evolution of  $\{R^* (\xb_T)\}_{T\ge 0}$ for $v$-processes on
        $\mathbb{S}^1$:
        $N=10^3$, $S_v(u)=u^{\alpha}$. Left: $\alpha =-1;$ $-0.53;$ $-0.05;$ $0.19;$
        $0.45;$ $0.66;$ $0.9,$ on the time interval $T\le 10 N.$ Right: $\alpha=1$, on the time interval $T\le 10^3 N.$}
    \label{Rstarevo}
\end{figure}

For a fixed configuration $\xb,$ let  $p_j := \lambda (C^{\xb}_{x_j})/\lambda(M),$ $
1\leq n\leq N$, be the probability that a random uniformly distributed point in $M$
lies in~$C_{x_{i}}^{\xb}$. Then Thiel's redundancy measure $R^*(\xb)$ of $\xb$ is
defined~\cite{LE79, OBS} as the entropy of the distribution $\{p_j\}_{j=1}^N$ relative
to the uniform one on $\{1,\dots, N\}$:
$$
R^*(\xb):= \text{ln}N+\sum_{j=1}^N p_j \ln p_j.
$$
Simulations show that this statistic works quite well over a reasonably large range of
$\alpha$ values: there is a very small variance in its values when $N\ge 10^3$, and a
strong enough dependence on~$\alpha$ to use $R^*$ for reliable estimation of the
parameter. Figure\,\ref{Rstar1D} displays (connected) average values of 25 independent
realizations for each of the values $\alpha = -2 + 0.3k,$ $0\le k\le 10,$ of
$R^*(\xb_T)$ after a large number of steps~$T$ (cf.~Fig.\,\ref{Rstarevo}).

The behaviour of the statistic $R^* (\xb_T)$ as $T$ increases could also be used to
estimate the rate of convergence to  stationarity. Figure\,\ref{Rstarevo} (left) gives
an indication of that rate for a range of $\alpha$ values. It shows that the rate is
quite high for $\alpha <0.5$, whereas when $\alpha$ approaches one, it takes the
process much longer to settle in a stable regime, and also that the oscillations of
$\{R^* (\xb_T)\}_{T\ge 0}$ are higher for those values of~$\alpha$. At the threshold
value $\alpha=1$, the sequence oscillates in quite diverse ranges of values
demonstrating metastable beahviour, see Fig.\,\ref{Rstarevo} (right).

\subsection{Two dimensions: edge effects}

The problem of `edge effects' (where the region boundary proximity can affect  the
shape/size of Voronoi cells) in simulations of Voronoi tessellations in (a rectangle)
$M\subset\R^2$ is commonly dealt with by treating the underlying space as a torus or by
considering only those parts of $M$ that lie within a window significantly distant from
the edge $\partial M$ of~$M$. This approaches is satisfactory for `static' point
processes~\cite{OBS}. In the case of the Voronoi processes, edge effects may propagate
in the course of the evolution, which presents an interesting problem by itself.

Let $\xb\subset M= [0,1]^2$ be finite, $A\subset \R^2$. The Voronoi tessellation
$\mathscr{T}_{\xb}$ induces a nearest-neighbour (NN) distance $d^{\xb}_{N\! N}(x_j,A)$
between $x_j\in \xb$ and $A$, defined as the length of the shortest `path' $x_j=x(1),
\dots, x(i)\in\xb$ with the properties that the Voronoi cells of $x(m)$ and $x(m+1)$
have a common edge, $1 \le m <i,$ and $C^{\xb}_{x(i)}\cap A\neq \varnothing$. We can
gauge the significance of edge effects by observing the changes in $\rstr$ when we
restrict attention to   $x_j\in\xb$ with $d^{\xb}_{N\!N}(x_j, \partial M)\ge m$ for a
fixed $m>0$. Edge effects were assessed by a two-way ANOVA using NN distance as one
factor, and the selection function as the other. In the case of the $v$-process with
$S_v(u)=u^{\alpha}$, we found that edge effects were only significant in the layer of
cells adjacent to the boundary, without any significant effect from~$\alpha$. For the
$n$-processes,   significant edge effects were detected  only for the `anti-few' and
`anti-many' selection functions (see subsection~\ref{s5} below), with no significant
propagation beyond the NN-depth of two.

To minimise edge effects in our results, the statistics were generally computed from the
set of cells with $d^{\xb}_{N\!N}(x_j,\partial M)\geq 3$.

\subsection{Two-dimensional $n$-processes}
\label{s5}

In this family of Voronoi processes, the evolution is driven by one of the most natural
local characteristics: the number of neighbours of the process points, \text{i.e.} the
number of the edges of their Voronoi cells. For non-degenerate random configurations in
$\mathbb R^2$, the vertices of the Voronoi tessellation \text{a.s.} terminate three
edges, so by Euler's theorem the mean number of edges for a cell is  six. There can be
cells with any number of edges $m\ge 3$, although in our simulations cells with $m> 13$
were rare, with relative frequency $<10^{-3}$.

We studied the $n$-processes on  $M=[0,1]^2$ with the following selection functions:
(i)~$S_n (n)=n$ (`vanilla', in which cells with a large number of neighbours are more
likely to be culled); (ii)~$S_n (n)=n^2$ (`anti-many',  which more severely penalises
cells with large numbers of neighbours); (iii)~$S_n (n)=(n-2)^{-2}$ (`anti-few', which
does just the opposite); (iv)~$S_n (n)=(0.1+|n-6|)^{-2}$ (`anti-6', which penalises
cells with six or close to six neighbours); (v)~$S_n (n)=|n-6|^2$ (`pro-6',which does
the opposite); and `sharp filters' focussing on cells with a given number neighbours:
(vi)~$S_n (5)=5000$ and $ S(n)=1$   if $n\neq 5$ (`anti-5'); and (vii)~$S_n (5)=1$ and
$S_n (n)=5000,$ if $n \neq 5$ (`pro-5'), and in addition the `anti-' and `pro-'
selectors for four and seven.

The statistics used in the study included:

\smallskip

(a)~the empirical probability mass function (EPMF) for the number of Voronoi NN's,

(b)~Thiel's redundancy measure, and also

(c)~the $J$-function of van Lieshout and Baddeley~\cite{BV95}.

\smallskip

The $J$-function $J(r)$ compares the `environment' of a `typical generator' of the
process with that of a `random point' in the underlying space, and is defined as the
ratio of the probabilities that a disk of radius $r$ centred at the given point is empty
of (other) generators of the process.

Formally, let $\xb$ be a stationary  isotropic point process in $\R^2$. Set
$$
 B(r):=\{x\in \mathbb R^2:\, |x|\le r\},\quad
 F(r):=\mathbf{P}\bigl(\xb\, \cap B(r)  \neq \varnothing \bigr),\quad
 G(r):= \mathbf{P}^{!0}\bigl({\xb} \cap B(r)  \neq \varnothing \bigr),
$$
where $\mathbf{P}^{!0}$ is the reduced Palm distribution of $\xb$ (the conditional
distribution of $\xb \backslash  \{0 \}$ given that there is a point at the origin).
Then the $J$-function is defined as
\begin{equation}
\label{jfunc}
 J(r):=\frac{1-G(r)}{1-F(r)} \quad \text{ for all } r\geq 0\text{ \ such that } F(r)<1.
\end{equation}
In the case of a Poisson point process of constant intensity, or complete spatial
randomness (CSR),  clearly $J(r)\equiv 1$, see e.g.~\cite{Amb02}. If there is a
tendency towards clustering obesrvable at scale $r$, then $J(r)<1$, while $J(r)>1$ if
there is a tendency toward regular spacing and $r$ is of the order of the average
distance between neighbours.

To estimate $J(r)$, one uses  empirical distribution functions $\hat{F}(r)$ and
$\hat{G}(r)$:
$$
\hat{J} (r) :=\frac{1-\hat{G}(r)}{1-\hat{F}(r)}.
$$
This is only reliably computable when $F(r) \leq 0.85$. Moreover, simulations showed
that  the values of~$\hat{J}$ have relatively large variances when $N\asymp 10^3$, so
we averaged the $\hat{J}$ values over 25--75  independent draws of the process,  and
performed a weighted cubic regression on the resulting points (with weights equal to
the reciprocals of the standard deviations; higher order regressions add little to the
goodness of fit).

\begin{figure}[!h]
    \centering
        \includegraphics{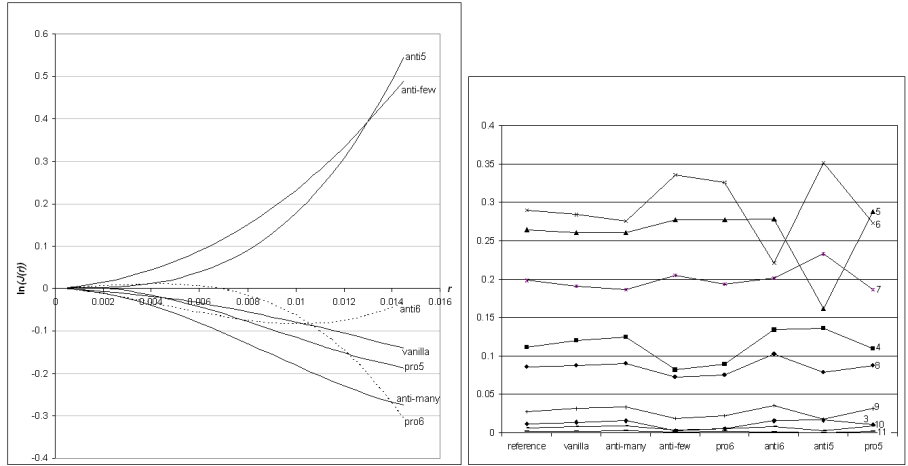}
        \caption{\footnotesize
        Comparison of $n$-processes ($N= 2000$, after $T=12N$ steps). Left: smoothed plots of  $\ln \hat{J} (r)$.
        Right: Proportions of cells with given numbers of neighbours (given
        at the right-hand end).
        Standard error is \text{approx.} 1\% of the given value in each case.}
    \label{nprocsComb}
\end{figure}

The results of our study for a selection of the above-listed $n$-processes are
summarised in Fig.\,\ref{nprocsComb} (statistics (c), (a)) and Table\,\ref{tabrst}
(statistic~(b)).

\begin{table}[!h]
\begin{center}
\begin{tabular}{c|c|c|c}
\hline
$\vphantom{\frac{1}{1}}$ $n$-process & $\overline R^*$ & s.d.($\rstr$) &s.e.($\overline R^*$)\\
\hline
vanilla& 0.1450& 0.0068& 0.0014\\
anti-many& 0.1548& 0.0073& 0.0015\\
anti-few& 0.1074& 0.0047& 0.0009\\
pro-6& 0.1615& 0.0083& 0.0012\\
anti-6& 0.1524& 0.0059& 0.0008\\
pro-5& 0.1479& 0.0067& 0.0013\\
anti-5& 0.1255& 0.0063& 0.0013
\end{tabular}
\end{center}
\caption{\footnotesize $\rstr$ values for various $n$-processes, with standard
deviations and standard errors. The values $\overline R^*$ are averages of 25
independent draws of $R^*$ ($N=2000$, after $T=12N$ steps). For CSR, $R^*=0.135.$}
\label{tabrst}
\end{table}

The $J$-function curves provide the most nuanced `picture' of the distribution
patterns, but these curves were obtained by averaging over a number of observations of
highly variable statistics and are not easily applied in practice. The statistic
$\rstr$ provides a crude but fairly sensitive measure of deviation from CSR. The NN
EPMF is probably the most effective statistic for the $n$-processes, but it doesn't
distinguish well between CSR and `anti-many' models (and also between `anti-few' and
`pro-6'), whereas the $J$-curves for these models are very different. The $J$-curves
for these and for the `sharp filters' on four and seven neighboured cells reveal the
interesting fact that selection functions favouring cells with $m>6$  ($m<6$) Voronoi
neighbours  produced more (less, resp.) uniform [than CSR] configurations.  The
`anti-few' and `anti-5' processes produce very similar $J$-curves, but are
distinguished by their effects on $\rstr$: `anti-few' gives $\rstr= 0.107$, while
`anti-5' has $\rstr=0.126$, which are significantly different from each other as well
as being significantly below the CSR value of 0.135. The low value for $\rstr$
corresponds to less variablity in cell areas, so this is consistent with the upward
$J$-curves. Notice also the unusual curvature of the `anti-6' curve, a result which was
confirmed by increasing the number of independent draws in this case to 75. In all
these cases, the $\rstr$ and $\hat J(r)$ statistics gave consistent results,
\text{i.e.} upward curves corresponded to $\rstr<0.135$.

\subsection{Two-dimensional $v$-processes}

\begin{figure}[!h]
    \centering
        \includegraphics{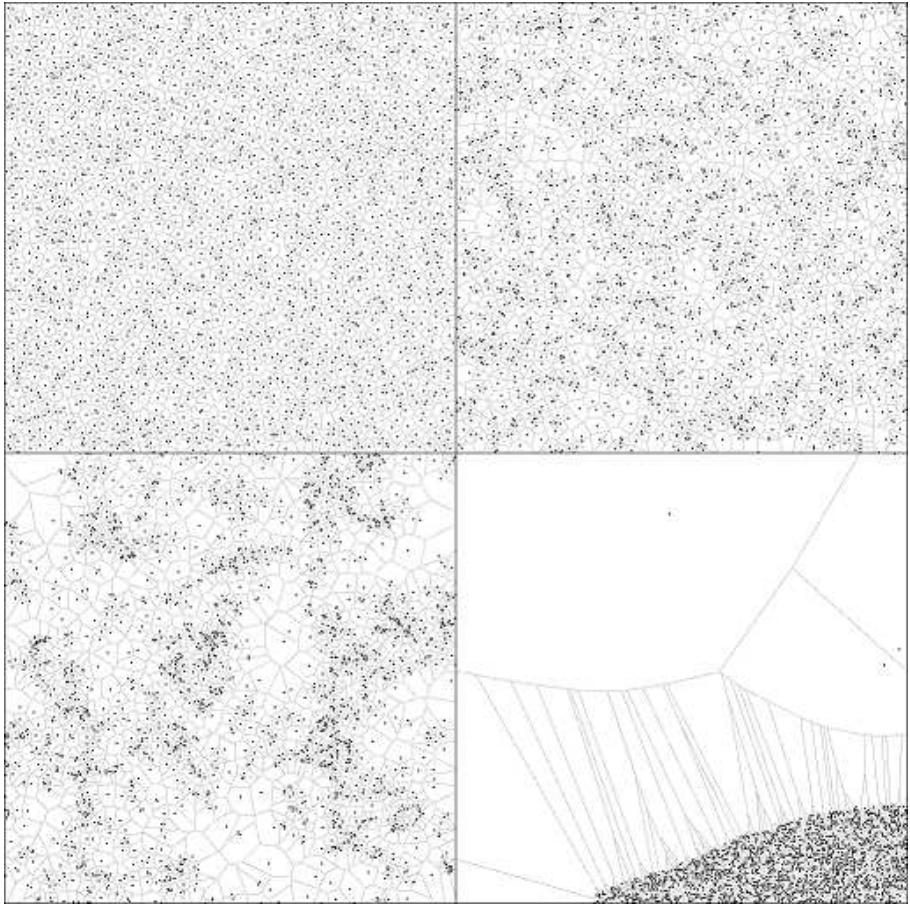}
        \caption{\footnotesize Realisations of the $v$-process with
        $S_v(u)=u^{\alpha}$ for some $\alpha$ values. Left-right, top-bottom:
        $\alpha=-3;$ 0.2; 1; 1.4. In each case, $N=2000$, $T=12N$.}
    \label{vprocmatrix}
\end{figure}

We study here the $v$-processes on $M=[0,1]^2$ with the scale-free selection functions
$S_v(u)=u^{\alpha}$. Figure\,\ref{vprocmatrix} shows a few examples of point patterns
after a considerable period of evolution. There is  a continuous gradation of patterns,
from almost regularly spaced ones for $\alpha<0$ through increasing levels of clustering
until we reach the threshold level at $\alpha=1$. For $\alpha>1$, the patterns become
increasingly concentrated in a small region, usually in a corner of the square, as the
number of steps increases.

\begin{figure}[!h]
    \centering
        \includegraphics{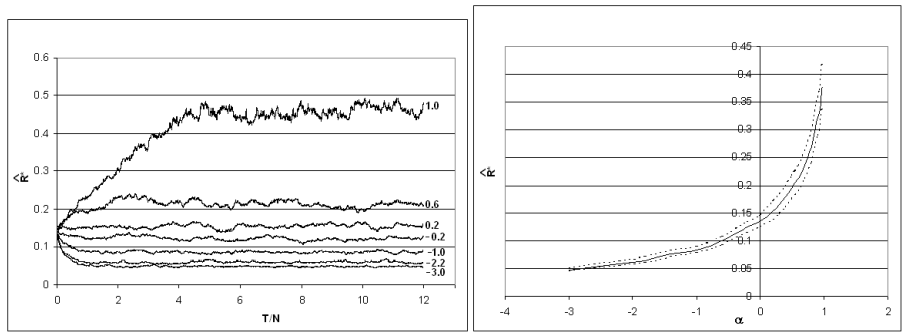}
        \caption{\footnotesize
        Thiel's redundancy measure for $v$-processes on $[0,1]^2$ ($N=2000$,
        $S_v(u)=u^{\alpha}$). Left: Evolution of $\{\rstr (x_T)\}_{T\le 12N}$ for
        $\alpha=-3;$ $-2.2$; $-1$;$-0.2$; $0.2$; $0.6$; 1.
        Right: $\rstr$ vs   $\alpha$ (after $T/N=12$ steps).}
    \label{vprocComb1}
\end{figure}

\begin{figure}[!h]
    \centering
        \includegraphics{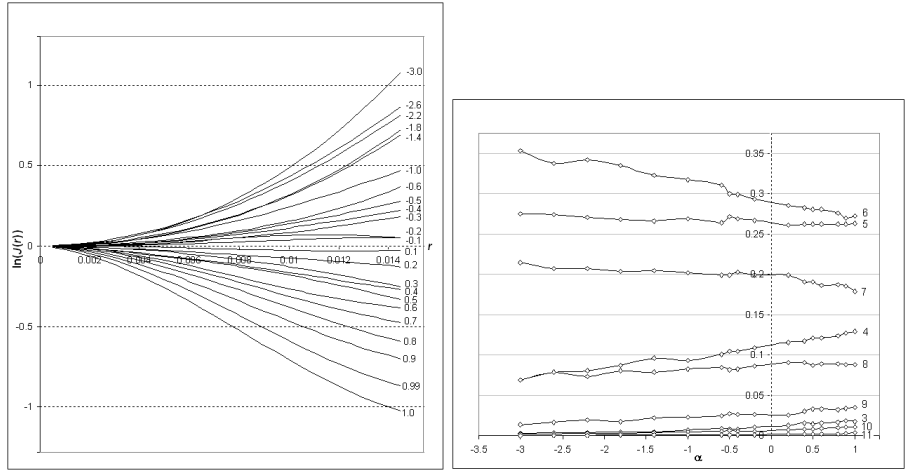}
        \caption{\footnotesize
        Graphical summary of results for $v$-processes on $[0,1]^2$ with
        $S_v(u)=u^{\alpha}$ for different $\alpha$ values ($N=2000$, $T=8N$).
        Left: Smoothed curves for   $\ln \hat J(r)$
        (numbers of independent draws range from 9 ($\alpha \leq -1.0$) to 59
        ($\alpha = 0.5, 0.6$)). Note the rapid change in the curve as $\alpha$ approaches 1.
        Right: Proportions of cells with given numbers of neighbours (given
        at the right-hand end) as
        functions of $\alpha$.}
    \label{vprocComb2}
\end{figure}

For $\alpha \leq 1$ we can employ the same statistics as for the $n$-process. As in the
one-dimensional case, Thiel's redundancy measure $\rstr$   serves well for inferring
the value of $\alpha$, and also provides a good indicator of the stabilisation of the
process, see~Fig.\,\ref{vprocComb1}. From Figs.\,\ref{vprocComb1},\,\ref{Rstarevo} it
appears  that the two-dimensional processes settle into equilibrium at roughly the same
rate as the one-dimensional processes. The $\ln \hat J(r)$ curves for $\alpha\in
(-3,1)$ form a fan as expected, see  Fig.\,\ref{vprocComb2}\,(left). More detailed
analysis reveals surprising (statistically significant) changes in curvature that occur
for $\alpha$ between $-0.3$ and $-0.2$ and between 0.5 and 0.6: it appears that a
subtle phase change does indeed occur around these values. In
Fig.\,\ref{vprocComb2}\,(right) we have plotted the NN EPMF as a function of $\alpha$.
Note that, as $\alpha$ increases (and so the amount of clustering), the proportions of
cells with different numbers of edges tend to come together; in other words, the cell
geometries become more variegated.

\section{Comparison with the area-interaction process}

The area-interaction point process (AIPP) of Baddeley and van Lieshout~\cite{BV95} is a
model which successfully produces a range of different clusterings. Kendall~\cite{K98}
presented a practical method for `perfect simulation' of this process. This is
summarised in Ambler~\cite{Amb02} from whose detailed account we derived our computer
program for simulations.

Let $M=[0,1]^2$ and $\mathscr{M}$ be the space of all finite configurations of points
in~$M$. The AIPP is specified by the density of its distribution with respect to that
of the unit rate Poisson process on $M$. Let $\lambda$ be Lebesgue measure on $M$,
$G:=B(\rho)$ for a fixed $\rho>0$, and  $\xb\oplus G := \left\{z:\, z=x+g, x\in \xb,
g\in G\right\}$  for $\xb\in \mathscr{M}$. The above density  is given by
$$
p(\xb)=C \beta^{\,\text{card} (\xb)}\gamma^{-\lambda(\xb\oplus G)},
$$
where $\beta, \gamma>0$ are parameters, $C$ a normalising constant. Note that the
Poisson process on $M$ with constant rate $\beta$ has density $\beta^{\,
\text{card}(\xb)}$. The parameter $\gamma$ defines the interactive component of the
process: when $0<\gamma <1$, configurations with high values of $\lambda(\xb\oplus G)$
for a given $\text{card}(\xb)$ are favoured, and so the interaction is described as {\em
repulsive\/},   while $\gamma >1$ favours $\xb$ with low $\lambda(\xb\oplus G)$ and so
constitutes the {\em attractive\/} case.

The expected number of points in $M$ is a complicated  function of $\gamma$, $\beta$
and $\rho$. We chose $\rho=0.01$, $\gamma=\gamma_{1}^{10^4}$ with $\gamma_{1}\in [0.3,
1.5]$, and  adjusted $\beta$ so that the number of points produced in the simulations
of the AIPP was close to 2000. This is because the AIPP doesn't re-scale in a simple
way, so comparison with the Voronoi $v$-process, in particular in regard to the
$J$-function, requires the samples to have roughly the same number of points.

\begin{figure}[h]
    \centering
        \includegraphics{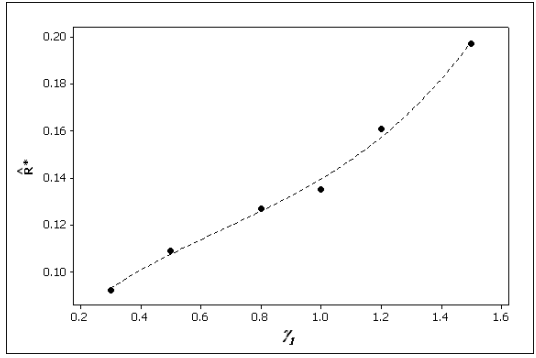}
        \caption{\footnotesize
        Plot of   $\rstr$ vs $\gamma_{1}$ for the AIPP
        as described above, with cubic regression line.
        Standard deviations of $\rstr$ are close to $0.005$. }
    \label{AIGR}
\end{figure}
\begin{figure}[h]
    \centering
        \includegraphics{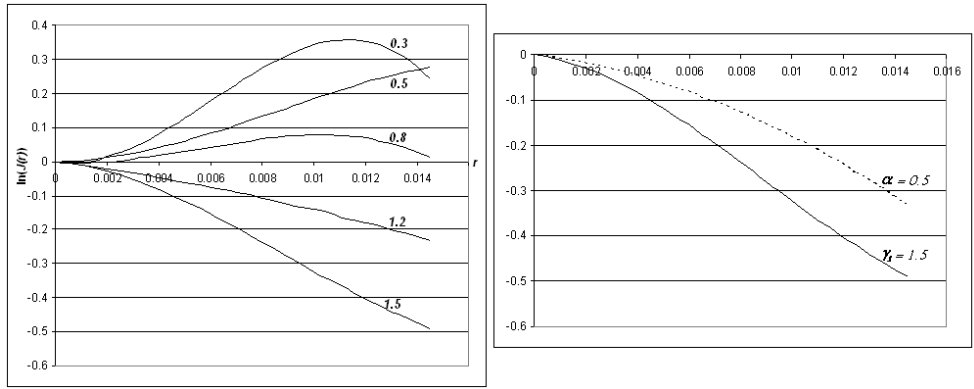}
        \caption{\footnotesize
        Left: Plots of smoothed $\ln \hat J(r) $ curves for the AIPP's
        with various values of $\gamma_{1}$, showing a range of attractive/repulsive effects.
        Right: The graph of $\ln \hat J(r)$ for the AIPP with $\gamma_1=1.5$ is compared to
        that for the
        scale-free Voronoi $v$-process with $\alpha=0.5$ (these two processes produce almost
        identical $\rstr$ values and NN EPMF's).}
    \label{AIJfns}
\end{figure}

\begin{figure}[!h]
    \centering
        \includegraphics{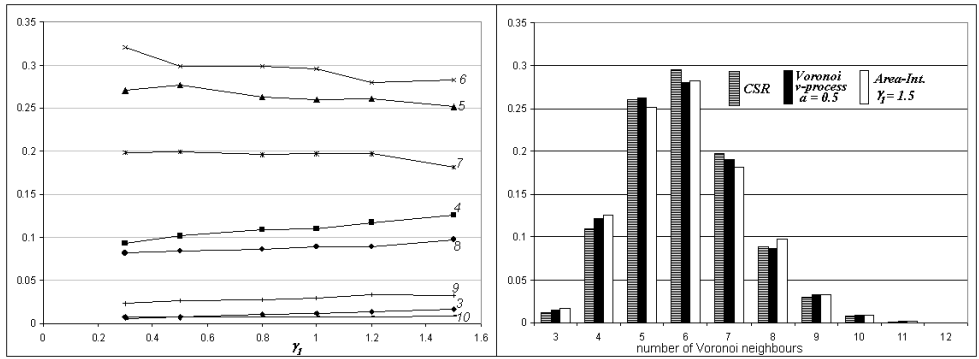}
        \caption{\footnotesize Left: Proportions of cells with given numbers of neighbours (given
        at the right-hand end) for  the AIPP's, plotted vs $\gamma_{1}$.
        Right: Comparison of the  NN EPMF's  for the AIPP with $\gamma_{1} = 1.5$, the
        Voronoi $v$-process with $\alpha = 0.5$,   and CSR.
        Standard errors are close to 0.005. }
    \label{AInn}
\end{figure}

Figures \ref{AIGR}--\ref{AInn} show the values of our three statistics for the AIPP. We
observe a continuous range of degrees of clustering on either side of CSR (when $\ln
J(r)\equiv 1$). It is interesting to compare these with those derived from the Voronoi
$v$-process. In particular, we can match results which produce close values for
$\rstr$. We have done this for the cases of the AIPP with $\gamma_{1} = 1.5$ and
Voronoi $v$-process with $\alpha = 0.5$ (for both processes, $\rstr\approx 0.2$). The
right-hand plots in Figs.\,\ref{AIJfns}, \ref{AInn} show that there is a significant
difference in the $J$-function curves, but that the Voronoi NN EPMF's do not yield
significant difference.

In summary, the $v$-processes and the AIPP's produce a range of clusterings depending on
a continuous parameter. The difference in the resulting point patterns can be detected
by the combination of our three statistics (a)--(c).


\end{document}